\documentclass[12pt,a4paper]{article}
\usepackage{mathrsfs}
\usepackage{amssymb}
\usepackage{amsmath}
\usepackage{amsfonts}
\usepackage{amstext}
\usepackage[top=1in, bottom=1in,left=0.75in,right=0.75in]{geometry}

\def \qed {\hfill \vrule height6pt width 6pt depth 0pt}

\begin{document}

\title{Eigenvalue, global bifurcation and positive solutions for a class of fully nonlinear problems
\thanks{Research supported by the NNSF of China (No. 11261052).}}
\author{{\small  Guowei Dai\thanks{Corresponding author. Tel: +86 931
7971124.\newline
\text{\quad\,\,\, E-mail address}: daiguowei@nwnu.edu.cn (G. Dai).%
} %\,\,\,\,Ruyun Ma
} \\
%EndAName
{\small Department of Mathematics, Northwest Normal
University, Lanzhou, 730070, PR China}\\
}
\date{}
\maketitle

\begin{abstract}
In this paper, we shall study global bifurcation phenomenon for the following Kirchhoff type problem
\begin{equation}
\left\{
\begin{array}{l}
-\left(a+b\int_\Omega \vert \nabla u\vert^2\,dx\right)\Delta u=\lambda u+h(x,u,\lambda)\,\,\text{in}\,\, \Omega,\\
u=0~~~~~~~~~~~~~~~~~~~~~~~~~~~~~~~~~~~~~~~~~~~~~~~~~\text{on}\,\,\Omega.
\end{array}
\right.\nonumber
\end{equation}
Under some natural hypotheses on $h$, we show that
$\left(a\lambda_1,0\right)$ is a bifurcation point of the above problem.
As applications of the above result, we shall determine the
interval of $\lambda$, in which there exist positive solutions for the above problem with $h(x,u;\lambda)=\lambda f(x,u)-\lambda u$,
where $f$ is asymptotically
linear at zero and is asymptotically 3-linear at infinity. To study global structure of bifurcation branch, we also
establish some properties of the first eigenvalue for a nonlocal eigenvalue problem. Moreover, we also provide a
positive answer to an open problem involving the case of $a=0$.
\end{abstract}\textbf{\ }

\numberwithin{equation}{section}

\numberwithin{equation}{section}

\section{Introduction}

\quad\, In this paper, we study global bifurcation phenomenon for the following  problem
\begin{equation}\label{bp}
\left\{
\begin{array}{l}
-\left(a+b\int_\Omega \vert \nabla u\vert^2\,dx\right)\Delta u=\lambda u+h(x,u,\lambda)\,\,\text{in}\,\, \Omega,\\
u=0~~~~~~~~~~~~~~~~~~~~~~~~~~~~~~~~~~~~~~~~~~~~~~~~~\text{on}\,\,\Omega,
\end{array}
\right.
\end{equation}
where $\Omega$ is a bounded domain in $\mathbb{R}^N$ with a smooth boundary $\partial\Omega$, $a>0$, $b>0$ are real constants, $\lambda$ is a parameter and
$h:\Omega\times \mathbb{R}^2\rightarrow\mathbb{R}$ satisfies the Carath\'{e}odory condition in the first two variable and
\begin{equation}\label{a1}
\lim_{ s\rightarrow0}\frac{h(x,s,\lambda)}{s}=0
\end{equation}
uniformly for a.e. $x\in\Omega$ and $\lambda$ on bounded sets. Moreover, we also assume that $h$ satisfies the growth restriction
\\

(G) There exist $c>0$ and $p\in\left(1,2^*\right)$ such that
\begin{equation}\label{a2}
\vert h(x,s,\lambda)\vert\leq c\left(1+\vert s\vert^{p-1}\right)\nonumber
\end{equation}
for a.e. $x\in \Omega$ and $\lambda$ on bounded sets, where
\begin{equation}\label{ef}
2^*=\left\{
\begin{array}{l}
\frac{2N}{N-2},~\,\, \text{if}\,\, N>2,\\
+\infty,~\,\, \text{if}\,\, N\leq 2.
\end{array}
\right.\nonumber
\end{equation}
\indent The problem (\ref{bp}) is related to the stationary problem of a model introduced by
Kirchhoff in 1883 to describe the transversal oscillations of a stretched string [\ref{K}]. More precisely, Kirchhoff proposed a model given by
the equation
\begin{equation}\label{kf}
\rho\frac{\partial^2u}{\partial
t^2}-\left(\frac{\rho_0}{h}+\frac{E}{2L}\int_0^L  \left\vert
\frac{\partial u}{\partial
x}\right\vert^2\,{d}x\right)\frac{\partial^2u}{\partial
x^2}=f(x,u),\nonumber
\end{equation}
where $\rho, \rho_0, h, E, L$ are constants, $f$ is the external force, which extends the
classical D'Alembert's wave equation, by considering the effect of
the changing in the length of the string during the vibration.
Problem (\ref{bp}) received much attention only after Lions [\ref{L}] proposed an
abstract framework to the problem. Some important and interesting
results can be found, for example, in [\ref{AP}, \ref{CCS}, \ref{DS}, \ref{DS1}, \ref{HZ}].
Recently, there are many mathematicians studying the
problem (\ref{bp}) by variational method, see [\ref{CKW}, \ref{CW}, \ref{MR}, \ref{MZ}, \ref{PZ}, \ref{ST}, \ref{ZP}] and the references therein.

Recently, the authors of [\ref{LLS}] studied problem (\ref{bp}) with $h(x,s,\lambda)=\lambda f(x,s)-\lambda s$ by using topological degree argument and variational method. Under some assumptions on $f$, they provided a positive answer to the existence of positive solutions to (\ref{bp}) for the cases $a>0$, $b>0$ and $a>0$, $b=0$. They pointed out that \emph{the situation of $a=0$ and $b>0$ is an open problem}.
The study of Kirchhoff type equations has already been extended to the case
involving the $p$-Laplacian and $p(x)$-Laplacian. We refer the readers to [\ref{APS}, \ref{CF}, \ref{DH}, \ref{DJ}, \ref{DL}, \ref{DM}, \ref{D}, \ref{D1}, \ref{F}] for an overview of and references on this subject.

A distinguishing feature of problem (\ref{bp}) is that the first equation
contains a nonlocal coefficient $a+b\int_\Omega \vert \nabla
u\vert^{2}\,dx$, and hence the equation is no longer a pointwise identity.
Moreover, the first equation of problem (\ref{bp}) with $h\equiv0$ even is not homogeneous. So problem (\ref{bp}) is a fully nonlinear problem which raises some essential difficulties to the study of this kind of problems.

The main aim of this paper is to establish the global bifurcation result for (\ref{bp}) and study its applications.
Let $\lambda_1$ denote the first eigenvalue of the following problem
\begin{equation}\label{le}
\left\{
\begin{array}{l}
-\Delta u=\lambda u\,\,\text{in}\,\, \Omega,\\
u=0~~~~~~~~\,\text{on}\,\,\Omega.
\end{array}
\right.
\end{equation}
It is well known that $\lambda_1$ is simple, isolated and the unique principal eigenvalue of (\ref{le}). The first main result of this paper is the following theorem.
\\ \\
\textbf{Theorem 1.1.} \emph{The pair $\left(a\lambda_1,0\right)$ is a bifurcation
point of (\ref{bp}). Moreover, there is a component $\mathcal{C}$ of the set of nontrivial solution of (\ref{bp}) in $\mathbb{R}\times H_0^1(\Omega)$
whose closure contains $\left(a\lambda_1, 0\right)$ and is either unbounded or contains a pair $\left(a\overline{\lambda}, 0\right)$
for some $\overline{\lambda}$, eigenvalue of (\ref{le}) with $\overline{\lambda}\neq\lambda_1$.}
\\ \\
\indent As is known to all, this result is proved in the case of $a=1$ and $b=0$ (see [\ref{R1}]). In fact, this result has been extended to the $p$-Laplacian problem
in the case of $a=1$ and $b=0$ (see [\ref{DPM}]). While, (\ref{bp}) is a full nonlinear problem.
So the Rabinowitz's global bifurcation theorem cannot be used directly to obtain our result. We shall transfer problem (\ref{bp}) into a new form and then use
Rabinowitz's global bifurcation theorem to prove Theorem 1.1 in Section 2

In order to find more detailed information of component $\mathcal{C}$ which is obtained in Theorem 1.1. In Section 3 we study the following eigenvalue problem
\begin{equation}\label{kee}
\left\{
\begin{array}{l}
-\left(\int_\Omega \vert \nabla u\vert^2\,dx\right)\Delta u=\mu u^3\,\, \text{in}\,\,\Omega,\\
u=0~~~~~~~~~~~~~~~~~~~~~~~~~~~\,\text{on}\,\,\partial\Omega.
\end{array}
\right.
\end{equation}
Let
\begin{equation}\label{fe}
\mu_1=\inf\left\{\left(\int_\Omega \vert \nabla u\vert^2\,dx\right)^2:u\in H_0^1(\Omega),\int_\Omega u^4\,dx=1\right\}.\nonumber
\end{equation}
We recall that a nodal domain of $u$ is a connected component
of $\Omega\setminus \{x\in \Omega:u(x)=0\}$.
Our second main theorem is:
\\ \\
\textbf{Theorem 1.2.} \emph{Assume that $N\leq 3$. Then $\mu_1$ is the first eigenvalue of (\ref{kee}) and has the following properties}\\

1. \emph{any eigenfunction $u$ corresponding to $\mu_1$ belongs to $C^{1,\alpha}\left(\overline{\Omega}\right)$ for some $\alpha \in(0,1)$, and $\frac{\partial u(x)}{\partial \gamma}<0$ if $u$ is nonnegative, where $\gamma$ is the outer unit normal at $x\in \partial \Omega$;}

2. \emph{$\mu_1$ is the principal eigenvalue of (\ref{kee});}

3. \emph{$\mu_1$ is simple in the sense that the eigenfunctions associated to it are merely a constant multiple of each other;}

4. \emph{(\ref{kee}) has a positive solution if and only if $\mu=\mu_1$;}

5. \emph{let $v$ be any eigenfunction associated to an eigenvalue $\mu>\mu_1$ and $\mathcal{N}$ be its any nodal domain. Then we have}
\begin{equation}\label{15}
\vert \mathcal{N}\vert\geq c,
\end{equation}
\emph{where $c>0$ is some constant depending only on $N$ and $\mu$;}

6. \emph{$\mu_1$ is isolated.}
\\

In [\ref{LLS}], the authors proved 3 and 4 in the case of $\Omega$ being a ball. However, we find that their proof contains a gap.
They claim that 2 has proved in [\ref{PZ}] which plays an essential role in their proof. This is not true. In fact, the authors of [\ref{PZ}] only proved the existence of $\mu_1$.
So the results of Theorem 1.2 not only fill the gap but also improve their results. To the best of our knowledge, the properties 1, 2, 5 and 6 are
the first results on this kind of problems.

In Section 4, we describe component $\mathcal{C}$ more detailed for problem (\ref{bp}) with $h(x,s,\lambda)=\lambda f(x,s)-\lambda s$, i.e., the following problem
\begin{equation}\label{fp}
\left\{
\begin{array}{l}
-\left(a+b\int_\Omega \vert \nabla u\vert^2\,dx\right)\Delta u=\lambda f(x,u)\,\,\text{in}\,\, \Omega,\\
u=0~~~~~~~~~~~~~~~~~~~~~~~~~~~~~~~~~~~~~~~~\,\text{on}\,\,\Omega.
\end{array}
\right.
\end{equation}
We assume that $f$ satisfies the following conditions
\\

(f1) $f:\Omega\times \mathbb{R}\setminus\mathbb{R}^-\rightarrow\mathbb{R}\setminus\mathbb{R}^-$ is a
Carath\'{e}odory function such that $f(x,s)s>0$ for a.e. $x\in \Omega$ and any $s>0$.

(f2) there exists $f_0\in (0,+\infty)$ such that
\begin{equation}
\lim_{s\rightarrow 0^+}\frac{f(x,s)}{a\lambda_1 s}=f_0\nonumber
\end{equation}
uniformly with respect to a.e. $x\in \Omega$.

(f3) there exists $f_\infty\in (0,+\infty)$ such that
\begin{equation}
\lim_{ s\rightarrow +\infty}\frac{f(x,s)}{b\mu_1 s^3}=f_\infty\nonumber
\end{equation}
uniformly with respect to a.e. $x\in \Omega$.\\

Our third main theorem is the following result.
\\ \\
\textbf{Theorem 1.3.} \emph{Assume that $N\leq 3$ and $f$ satisfies (f1)--(f3). Then $\left(1/f_0, 0\right)$ is the unique bifurcation point of the set of positive solution of (\ref{fp}). Moreover, there is an unbounded component $\mathcal{C}$ in $\mathbb{R}\times H_0^1(\Omega)$
whose closure links $\left(1/f_0, 0\right)$ to $\left(1/f_\infty, \infty\right)$.}
\\ \\
\textbf{Remark 1.1.} The reason of needing $N\leq 3$ in Theorem 1.2 is to insure $u^3\in L^2(\Omega)$ for any $u\in H_0^1(\Omega)$ which is
crucial for our proof. However, we doubt its necessity for the theorem. In addition, we pose the condition $N\leq 3$ because our proof depends Theorem 1.2 and the fact of $4<2^*$.
We also doubt its necessity for Theorem 1.3.\\
\\
\textbf{Remark 1.2.} Note that Theorem 1.1 of [\ref{LLS}] is just our corollary of Theorem 1.3. We also note
that Theorem 1.2 (ii) (with $a>0$) of [\ref{LLS}] is the corollary of Theorem 1.3. Clearly, Theorem 1.2 (i) of [\ref{LLS}] does not occur because
$\left(1/f_\infty, \infty\right)$ is the unique bifurcation point of positive solution set of (\ref{fp}) where a bifurcation from infinity occurs.
So Theorem 5.1 and 5.2 of [\ref{LLS}] are also our corollary of Theorem 1.3.
Moreover, our assumptions on $f$ and $\Omega$ are more concise and weaker than the corresponding
ones of [\ref{LLS}, Theorem 1.1, 1.2, 5.1 and 5.2].
\\
%
%Furthermore, we also have the following theorem.
%\\ \\
%\textbf{Theorem 1.4.} \emph{Assume that $a>0$ and $N\leq 3$. Besides the assumptions (A1)--(A3), we also assume that $f$ satisfies}\\
%
%(f4) \emph{$f(x,\cdot)\in C^1\left(\mathbb{R}^+,\mathbb{R}\setminus\mathbb{R}^+\right)$ such that $f'(x,s)s<f(x,s)$ for any $x\in \Omega$ and $s>0$.}\\
%
%\noindent \emph{Then (\ref{fp}) has a unique positive solutions if and only if $\lambda\in\left(\frac{1}{f_0},\frac{1}{f_\infty}\right)\cup\left(\frac{1}{f_\infty},\frac{1}{f_0}\right)$. Moreover, all positive solutions $u(\lambda,\cdot)$ lie on a smooth curve $\mathcal{C}$
%whose closure links $\left(1/f_0, 0\right)$ to $\left(1/f_\infty, \infty\right)$.}
%\\
\\
\textbf{Remark 1.3.} From Theorem 1.3, we can see that any positive of (\ref{fp}) lies in $\mathcal{C}$. That is to say,
we find the range of all positive solutions. So it only needs to study the structure and formulation of $\mathcal{C}$ to find the positive of (\ref{fp}).
\\ \\
\textbf{Remark 1.4.} By standard elliptic regularity theory (see [\ref{F1}, \ref{FZ}]), we know that any weak solution of (\ref{bp}) or (\ref{fp}) belongs to $C^{1,\alpha}(\overline{\Omega})$ with $\alpha\in(0,1)$ under the condition of (\ref{a1}) and (G) or (f1)--(f3).\\

In Section 5, we consider the case of $a=0$ in (\ref{fp}) and give a positive answer to the above mentioned open problem.
Concretely, we give an assumption on $f$ as follows \\

(f4) there exists $\widetilde{f}_0\in (0,+\infty)$ such that
\begin{equation}
\lim_{s\rightarrow 0^+}\frac{f(x,s)}{\lambda_1 s}=\widetilde{f}_0\nonumber
\end{equation}
uniformly with respect to a.e. $x\in \Omega$.\\

Our last main result is the following theorem.
\\ \\
\textbf{Theorem 1.4.} \emph{Assume that $a=0$, $N\leq 3$ and $f$ satisfies (f1), (f3) and (f4). Then $\left(0, 0\right)$ is the unique bifurcation point of the set of positive solution of (\ref{fp}). Moreover, there is an unbounded component $\mathcal{C}$ in $\mathbb{R}\times H_0^1(\Omega)$
whose closure links $\left(0, 0\right)$ to $\left(1/f_\infty, \infty\right)$.}\\
\\
\textbf{Remark 1.5.} From Theorem 1.4, we can easily see that (\ref{fp}) has a positive solution with $\lambda=1$ and $f_\infty<1$.
This give a positive answer to an open problem proposed by the authors of [\ref{LLS}].
\\ \\
\textbf{Remark 1.6.} Note that $f$ is asymptotically linear at zero in Theorem 1.4 which is different from [\ref{LLS}].
In [\ref{LLS}], the $f$ is assumed to be asymptotically 3-linear at zero in the case of $a=0$. We shall consider this situation in our future work.
\\

The last Section concludes the paper and outlines our future work. We end this section by introducing some notation conventions which will be used later in this paper. Let $X$ be the usual Sobolev space $H_0^1(\Omega)$ with the norm $\Vert u\Vert=\left(\int_\Omega \vert \nabla u\vert^2\,dx\right)^{1/2}$ and $X^*$ be its dual space. Denote by $\langle\cdot,\cdot\rangle$ the duality pairing between $X$ and $X^*$.
We write $u_n\rightharpoonup u$ and $u_n\rightarrow u$ the weak convergence
and strong convergence of sequence $\left\{u_n\right\}$ in $X$, respectively. Use $q'=q/(q-1)$ to denote the conjugative number of $q$ with $q>1$. For a measurable set $A$ of $\mathbb{R}^N$, we denote its measure by $\vert A\vert$. Also, denote by $c$ and $c_{i}$, $i\in \mathbb{N}$, the generic positive constants (the exact value may be different from line to line).

\section{Global bifurcation}

\quad\, Firstly, consider the following auxiliary problem
\begin{equation}\label{ap}
\left\{
\begin{array}{l}
-\Delta u=f(x)\,\, \text{in}\,\,\Omega,\\
u=0~~~~~~~~~~\text{on}\,\,\partial\Omega.
\end{array}
\right.
\end{equation}
As is known to all, problem (\ref{ap}) possesses a unique weak solution for each $f\in X^*$.
Let us denote by $G(f)$ the unique weak solution of (\ref{ap}). Then $G:X^*\rightarrow X$ is a
linear continuous operator. Since the embedding of $X\hookrightarrow L^q(\Omega)$ is compact for each $q\in\left(1,2^*\right)$,
the restriction of $G$ to $L^{q'}(\Omega)$ is a completely continuous operator.

Clearly, the pair $(\lambda,u)$ is a solution of (\ref{bp}) if and only if $(\lambda,u)$ satisfies
\begin{equation}\label{oe}
u=G\left(\frac{1}{a+b\Vert u\Vert^2}(\lambda u+H(\lambda,u))\right),
\end{equation}
where $H(\lambda,\cdot)$ denotes the usual Nemitsky operator associated with $h$.
From condition (G) and noting $2<2^*$, we can see that $G:L^{p'}(\Omega)\cup L^{2}(\Omega)\rightarrow X$ is completely continuous.
\\ \\
\textbf{Proof of Theorem 1.1.} Let
\begin{equation}
Lu=\frac{1}{a}G(u),\,\, \widetilde{H}(\lambda,u)=\frac{1}{a+b\Vert u\Vert^2} G(H(\lambda,u))-\frac{\lambda b\Vert u\Vert^2}{a\left(a+b\Vert u\Vert^2\right)}G(u).\nonumber
\end{equation}
Clearly, $L:X\rightarrow X$ is linear completely continuous, $\widetilde{H}:\mathbb{R}\times X\rightarrow X$ is compact.
Moreover, it is easy to see that $a\lambda_1$ is simple characteristic value of $L$. Then equation (\ref{oe}) is equivalent to
\begin{equation}
u=\lambda Lu+\widetilde{H}(\lambda,u).\nonumber
\end{equation}
Next, we show that $\widetilde{H}=o(\Vert u\Vert)$ at $u=0$ uniformly on bounded $\lambda$ intervals.
It is sufficient to show that
\begin{equation}
\lim_{\Vert u\Vert\rightarrow 0}\frac{H\left(x,u\right)}{\left\Vert u\right\Vert}=0\,\, \text{in} \,\, L^{p'}(\Omega).\nonumber
\end{equation}
Without loss of generality, we may assume that $p>2$. Otherwise, we can consider $\tilde{p}=cp$, $c>1$ such that $\tilde{p}\in\left(2,2^*\right)$.
From $p<2^*$, we can see that
\begin{equation}
\frac{p'(p-2)}{2^*}<\frac{2^*-p'}{2^*}.\nonumber
\end{equation}
So we can choose a real number $r>1$ such that
\begin{equation}
\frac{p'(p-2)}{2^*}\leq\frac{1}{r}\leq\frac{2^*-p'}{2^*}.\nonumber
\end{equation}
It follows that
\begin{equation}\label{p2r}
p' r(p-2)\leq 2^*\,\,\text{and}\,\, p' r'\leq2^*.
\end{equation}
\indent For any $\varepsilon>0$, in view of (\ref{a1}) and (G), we can choose positive numbers $\delta=\delta(\varepsilon)$ and $M=M(\delta)$ such that
for a.e. $x\in \Omega$, the following relations hold:
\begin{equation}
\left\vert \frac{h(x,s,\lambda)}{s}\right\vert\leq \varepsilon\,\,\,\text{for}\,\, 0<\vert s\vert\leq \delta,\nonumber
\end{equation}
\begin{equation}
\left\vert \frac{h(x,s,\lambda)}{s}\right\vert\leq M\vert s\vert^{p-2}\,\,\,\text{for}\,\, \vert s\vert>\delta.\nonumber
\end{equation}
Then we can obtain that
\begin{equation}
\int_\Omega\left\vert \frac{H\left(\lambda,u\right)}{u}\right\vert^{p'r}\,dx\leq \varepsilon \vert \Omega\vert+M^{p' r}\int_\Omega \left\vert u\right\vert^{p'r(p-2)}\,dx.\nonumber
\end{equation}
From this inequality, (\ref{p2r}) and $u\rightarrow 0$ in $X$, we get that
\begin{equation}\label{27}
\left\vert \frac{H\left(\lambda,u\right)}{u}\right\vert^{p'}\rightarrow 0\,\,\text{in}\,\, L^r(\Omega).
\end{equation}
Let $v=u/\Vert u\Vert$. By the boundedness of $v$ in $X$, (\ref{p2r}) and the continuous embedding of $X\hookrightarrow L^{2^*}(\Omega)$, we have that
\begin{equation}\label{28}
\int_\Omega \left\vert v\right\vert^{p'r'}\,dx\leq c
\end{equation}
for some constant $c>0$.
Then from (\ref{27}), (\ref{28}) and the H\"{o}lder's inequality, we obtain that
\begin{eqnarray}
\int_\Omega \left\vert \frac{H\left(\lambda,u\right)}{\left\Vert u\right\Vert}\right\vert^{p'}\,dx&=&\int_\Omega \left\vert \frac{H\left(\lambda,u\right)}
{\left\vert u\right\vert}\right\vert^{p'} \left\vert v\right\vert^{p'}\,dx\nonumber\\
&\leq &\left(\int_\Omega\left\vert \frac{H\left(\lambda,u\right)}{u}\right\vert^{p'r}\,dx\right)^{1/r}\left(\int_\Omega \left\vert v\right\vert^{p'r'}\,dx\right)^{1/r'}\nonumber\\
&\rightarrow&0.\nonumber
\end{eqnarray}
Now, from global bifurcation theory (see [\ref{R2}, Theorem 1.3]), we get the existence of
a global branch of the set of nontrivial solution of (\ref{bp}) emanating from
$\left(a\lambda_1, 0\right)$.\qed

\section{Properties of the first eigenvalue of a nonlocal problem}

\quad\, In order to study more detailed information of the component $\mathcal{C}$ which is obtained in Theorem 1.1, we must consider the
following eigenvalue problem
\begin{equation}\label{ke}
\left\{
\begin{array}{l}
-\left(\int_\Omega \vert \nabla u\vert^2\,dx\right)\Delta u=\mu u^3\,\, \text{in}\,\,\Omega,\\
u=0~~~~~~~~~~~~~~~~~~~~~~~~~~~\,\text{on}\,\,\partial\Omega.
\end{array}
\right.
\end{equation}
Let
\begin{equation}
I(u)=\Vert u\Vert^4,\,\,u\in S:=\left\{u\in X:\int_\Omega u^4\,dx=1\right\}.\nonumber
\end{equation}
Denote by $\mathcal{A}$ the class of closed symmetric subsets of $S$, let
\begin{equation}
\mathcal{F}_m=\left\{A\in \mathcal{A}:i(A)\geq m-1\right\},\nonumber
\end{equation}
where $m$ is a positive inter number and $i(A)$ denotes the Yang index of $A$.
The authors of [\ref{PZ}] have proved that problem (\ref{ke}) possesses
an unbounded sequences of minimax eigenvalues
\begin{equation}
0<\mu_1\leq \mu_2<\cdots\nonumber
\end{equation}
such that
\begin{equation}
\mu_m=\inf_{A\in \mathcal{F}_m}\max_{u\in A} I(u).\nonumber
\end{equation}
In particular, if $m=1$, taking $A=\{u,-u: u\in S\}$, we can get that
\begin{equation}\label{fe}
\mu_1=\inf\left\{\Vert u\Vert^4:u\in X,\int_\Omega u^4\,dx=1\right\}.\nonumber
\end{equation}
\indent Next, we are going to study the properties of $\mu_1$.
These properties are important in the studying of the global bifurcation phenomena.
\\ \\
\textbf{Proposition 3.1.} \emph{Let $u$ be any eigenfunction corresponding to $\mu_1$. Then $u\in C^{1,\alpha}\left(\overline{\Omega}\right)$ for some $\alpha \in(0,1)$. Furthermore, $\frac{\partial u(x)}{\partial \gamma}<0$ if $u$ is nonnegative, where $\gamma$ is the outer unit normal at $x\in \partial \Omega$.}
\\
\\
\textbf{Proof.} Note that problem (\ref{ke}) is homogeneous. So by scaling we may suppose that $\Vert u\Vert=1$. It follows that
\begin{equation}\label{20131214}
\left\{
\begin{array}{l}
-\Delta u=\mu_1 u^3\,\, \text{in}\,\,\Omega,\\
u=0~~~~~~~~~~\,\text{on}\,\,\partial\Omega.
\end{array}
\right.
\end{equation}
By the embedding of $X\hookrightarrow L^{2^*}$ and $N\leq 3$, we can have $\mu_1 u^3(x):=f(x)\in L^2(\Omega)$. By [\ref{GT}, Theorem 8.12], we know that
$u\in W^{2,2}(\Omega)$. Furthermore, by the general Sobolev embedding theorem [\ref{E}, p. 270], we get $u\in C^\gamma\left(\overline{\Omega}\right)$ for some $\gamma \in(0,1)$.
Moreover, by the definition of weak derivative, we can get that $\nabla f= \mu_13 u^2 \nabla u$. It is easy to verify that $\nabla f \in L^2(\Omega)$. So we have $f\in W^{1,2}(\Omega)$.
By [\ref{GT}, Theorem 8.13], we know that
$u\in W^{3,2}(\Omega)$. Again using the general Sobolev embedding theorem [\ref{E}, p. 270], we get $u\in C^{1,\alpha}\left(\overline{\Omega}\right)$ for some $\alpha \in(0,1)$. Furthermore, if $u\geq0$, by the strong maximum principle [\ref{FZZ}, Theorem 1.2], we get $\frac{\partial u(x)}{\partial \gamma}<0$ for all $x\in \partial \Omega$.
\qed
\\ \\
\textbf{Proposition 3.2.} \emph{Let $u$ be an eigenfunction associated with $\mu_1$, then either $u>0$ or $u<0$ in $\Omega$, i.e., $\mu_1$ is the principal eigenvalue of (\ref{ke}).}
\\ \\
\textbf{Proof.} We notice that if $u$ is an eigenfunction, so is $v:=\vert u\vert$. Without loss of generality, we assume that $\Vert v\Vert=1$.
So we have
\begin{equation}
\left\{
\begin{array}{l}
-\Delta v=\mu_1 v^3\,\, \text{in}\,\,\Omega,\\
v=0~~~~~~~~~~\,\text{on}\,\,\partial\Omega.
\end{array}
\right.\nonumber
\end{equation}
By the strong maximum principle [\ref{GT}, Theorem 8.19], we know that $v>0$ in the whole domain. By the continuity of $u$, either $u$ or $-u$ is
positive in the whole domain.\qed
\\ \\
\textbf{Proposition 3.3.} \emph{The principal eigenvalue $\mu_1$ is simple, i.e., if $u$ and $v$ are two eigenfunctions
associated with $\mu_1$, then there exists $c$ such that $u = cv$.}
\\ \\
\textbf{Proof.} Define $f^1$ and $f^2$ on $X$ by
\begin{equation}
f^1(w)=\frac{1}{4}\Vert w\Vert^4,\,\, f^2(w)=\frac{1}{4}\int_\Omega w^4\,dx.\nonumber
\end{equation}
Set $J_\mu(w)=f^1(w)-\mu f^2(w)$ and $Q(w)=f^2(w)/f^1(w)$, then
\begin{equation}
J_\mu(w)\underset{=}{<}0\,\,\text{if and only if}\,\, Q(w)\underset{=}{>}\frac{1}{\mu}.\nonumber
\end{equation}
It follows that there exists $u_0\in X$ such that $u_0\not\equiv 0$ and
\begin{equation}
\min\left\{J_{\mu_1}(w):w\in X,w\not\equiv 0\right\}=0=J_{\mu_1}\left(u_0\right).\nonumber
\end{equation}
Consequently, $u$ becomes a nontrivial solution of (\ref{ke}) with $\mu=\mu_1$. Conversely, if $u$ is a solution of (\ref{ke}) with $\mu=\mu_1$,
then multiplication of (\ref{ke}) by $u$ gives $J_{\mu_1}(u)=0$. Thus we find that
\begin{equation}\label{iff}
u\,\,\text{is a solution of (\ref{ke}) with}\,\, \mu=\mu_1\,\,\text{if and only if}\,\, J_{\mu_1}(u)=0.
\end{equation}
\indent Let $u$, $v$ be two eigenfunctions
associated with $\mu_1$ and put
\begin{equation}
M(t,x)=\max\{u(x),tv(x)\}\,\, \text{and}\,\, m(t,x)=\min\{u(x),tv(x)\}\nonumber
\end{equation}
for all $t>0$.
Without loss of generality, we assume that $\Vert u\Vert=\Vert v\Vert=1$.
Set
\begin{equation}
(u\vee v)(x)=\max\{u(x),v(x)\}\,\, \text{and}\,\, (u\wedge v)(x)=\min\{u(x),v(x)\},\nonumber
\end{equation}
and let
\begin{equation}
\Omega_1=\{x\in\Omega:(u\vee v)(x)=u(x)\}\,\, \text{and}\,\, \Omega_2=\{x\in\Omega:(u\vee v)(x)=v(x)\}.\nonumber
\end{equation}
Clearly, we have $\Omega=\Omega_1\cup \Omega_2$.

Define
\begin{equation}
g^1(w)=\Vert w\Vert^2,\,\, g^2(w)=\int_\Omega \left(\frac{w}{\sqrt{\Vert w\Vert}}\right)^4\,dx.\nonumber
\end{equation}
By some simple computation, we can obtain that
\begin{equation}\label{g1}
g^1(u\vee v)+g^1(u\wedge v)\leq g^1(u)+g^1(v)
\end{equation}
and
\begin{equation}\label{g2}
g^2(u\vee v)+g^2(u\wedge v)\geq g^2(u)+g^2(v).
\end{equation}
\indent Define $\Phi_\mu(w)=g^1(w)-\mu g^2(w)$. It is easy to see that $J_{\mu_1}(w)=\frac{1}{4}\Vert w\Vert^2\Phi_{\mu_1}(w)\geq 0$.
Hence, by (\ref{g1}) and (\ref{g2}), we get that
\begin{equation}
0\leq \Phi_{\mu_1}(M)+\Phi_{\mu_1}(m)\leq\Phi_{\mu_1}(u)+\Phi_{\mu_1}(tv)=4\left(J_{\mu_1}(u)+\frac{J_{\mu_1}(tv)}{t^2}\right)=4\left(J_{\mu_1}(u)+t^2J_{\mu_1}(v)\right)=0.\nonumber
\end{equation}
So we have that
\begin{equation}
\Phi_{\mu_1}(M)=\Phi_{\mu_1}(m)=0.\nonumber
\end{equation}
It is not difficult to show that $\Vert M(t,\cdot)\Vert\leq 1+t^2$ and $\Vert m(t,\cdot)\Vert\leq 1+t^2$. So we get
\begin{equation}
J_{\mu_1}(M)=J_{\mu_1}(m)=0.\nonumber
\end{equation}
Hence, by (\ref{iff}), $M$ and $m$ turn out to be solutions of (\ref{ke}) with $\mu=\mu_1$ for all $t>0$.

We note that since $u$, $v\in C^{1,\alpha}\left(\overline{\Omega}\right)$ (comes from Proposition 3.1), $u$, $v$ are absolutely continuous in each
variable (on segments in $\Omega$) for almost all values of other variables, and their partial and generalized derivatives coincide almost everywhere, (see [\ref{M}]).
By Proposition 3.2, we know that $u\neq 0$ and $v\neq 0$. Without loss of generality, we assume that $u$ and $v$ are positive in $\Omega$.
Moreover, by virtue of the facts that $v>0$ in $\Omega$ and $u$, $v\in C^{1,\alpha}\left(\overline{\Omega}\right)$, we have $u/v$ belongs to $C(\Omega)\cap W_{loc}^{1,2}(\Omega)$.

For a.e. $x_0\in \Omega$, we set $t_0=u\left(x_0\right)/v\left(x_0\right)>0$. Then, for any unit vector $e$, we get that
\begin{equation}
u\left(x_0+he\right)-u\left(x_0\right)\leq M\left(t_0,x_0+he\right)-M\left(t_0,x_0\right),\nonumber
\end{equation}
\begin{equation}
t_0v\left(x_0+he\right)-t_0v\left(x_0\right)\leq M\left(t_0,x_0+he\right)-M\left(t_0,x_0\right).\nonumber
\end{equation}
Dividing above two inequalities by $h>0$ and $h<0$ and letting $h$ tend to $\pm 0$, we get
\begin{equation}
\nabla_x u\left(x_0\right)=\nabla_x M\left(t_0,x_0\right)=t_0\nabla_x v\left(x_0\right).\nonumber
\end{equation}
Thus
\begin{equation}
\nabla_x\left(\frac{u}{v}\left(x_0\right)\right)=\frac{\nabla_x u\left(x_0\right)v\left(x_0\right)-u\left(x_0\right)\nabla_x v\left(x_0\right)}{v^2\left(x_0\right)}=0.\nonumber
\end{equation}
Therefore, we obtain that $u(x)/v(x)=c$ in $\Omega$ for some constant $c\neq 0$.\qed
\\ \\
\textbf{Proof of Theorem 1.2.} By Proposition 3.1--3.3, we only need to prove property 4, 5 and 6.

4. Suppose on the contrary that (\ref{ke}) with $\mu>\mu_1$ has a positive solution $v$, and let $u$ be a positive eigenfunction
corresponding to $\mu_1$. Similarly to Proposition 3.1, we can show that
$v\in C^{1,\alpha}\left(\overline{\Omega}\right)$ for some $\alpha \in(0,1)$ and $\frac{\partial v(x)}{\partial \gamma}<0$. By this, Proposition 3.1 and the fact that $t v$ is also an eigenfunction of (\ref{ke}), we may assume that $u\leq v$.
Let $A=\left(f^1\right)'$ and $B=\left(f^2\right)'$, where $f^1$ and $f^2$ come from the proof of Proposition 3.3.
Then $w$ is the weak solution of (\ref{ke}) if and only if
\begin{equation}
Aw=\mu Bw.\nonumber
\end{equation}
It is not difficult to show that $Bu\leq Bv$.
Then we get that
\begin{equation}
Au=\mu_1Bu\leq \mu_1Bv=\mu B(\eta v)=A(\eta v)\,\, \text{with} \,\, \eta=\left(\frac{\mu_1}{\mu}\right)^{1/3}<1.\nonumber
\end{equation}
It follows from Proposition 5.1 which will be proved in section 5 that $u\leq \eta v$. Repeating this argument $n$ times, we gain that $0\leq u\leq \eta^nv$, Letting $n\rightarrow +\infty$, we get
$u\equiv 0$. This is a contradiction. So $v$ must change sign.

5. By an argument similar to that of Proposition 3.1, we know that $v\in C(\Omega)$ then $v_{|_{\mathcal{N}}}\in H_0^1(\mathcal{N})$. Define
\begin{equation}
w=\left\{
\begin{array}{l}
v, \,\,x\in \mathcal{N},\\
0,\,\,x\in \Omega\setminus\mathcal{N}.
\end{array}
\right.\nonumber
\end{equation}
It is easy to see $w\in X$. We first consider the case of $N=3$. Then we have that
\begin{equation}
\left(\int_\mathcal{N} \left\vert \nabla w\right\vert^2\,dx\right)^2=\mu \int_\mathcal{N}w^4\,dx.\nonumber
\end{equation}
By the H\"{o}lder inequality and the Sobolev embeddings we have that
\begin{equation}
c^4\left(\int_\mathcal{N}\vert w\vert^6\,dx\right)^{2/3}\leq\left(\int_\mathcal{N} \left\vert \nabla w\right\vert^2\,dx\right)^2=\mu \int_\mathcal{N} w^4\,dx
\leq\mu\left(\int_\mathcal{N}w^6\,dx\right)^{2/3} \vert\mathcal{N}\vert^{1/3},\nonumber
\end{equation}
where $c>0$ is the best embedding constant of $H_0^1\left(\mathcal{N}\right)\hookrightarrow L^6\left(\mathcal{N}\right)$.
It follows that
\begin{equation}
\vert\mathcal{N}\vert\geq \frac{c^{12}}{\mu^3}.\nonumber
\end{equation}
For $N=1$ or $2$, it is well known that $H_0^1\left(\mathcal{N}\right)$ is continuously embedding $L^\infty\left(\mathcal{N}\right)$ with the best embedding constant $c>0$.
From this fact and reasoning as above, we can show that
\begin{equation}
\vert\mathcal{N}\vert\geq \frac{c^4}{\mu}.\nonumber
\end{equation}
\indent 6. Assume by contradiction that there exists a sequence of eigenvalues
$\mu_n\in\left(\mu_1, \delta\right)$ for some constant $\delta>\mu_1$ which converges to $\mu_1$. Let $u_n$ be the corresponding eigenfunctions.
Property 4 implies that $u_n$ changes sign.
Integration by parts helps to yield
\begin{equation}
\left(\int_\Omega\left\vert \nabla u_n\right\vert^2\,dx\right)^2=\mu_n\int_\Omega \left( u_n\right)^4\,dx.\nonumber
\end{equation}
Define
\begin{equation}
v_n:=\frac{u_n}{\left(\int_\Omega \left( u_n\right)^4\,dx\right)^{1/4}}.\nonumber
\end{equation}
Obviously, $v_n$ is bounded in $X$ so there exists a subsequence, denoted again by $v_n$, and $v\in X$ such that
$v_n\rightharpoonup v$ in $X$ and $v_n\rightarrow v$ in $L^4(\Omega)$.
Since functional $I$ is sequentially weakly lower semi-continuous, we have that
\begin{equation}
\left(\int_\Omega \left\vert \nabla v\right\vert^2\,dx\right)^2 \leq\liminf_{n\rightarrow+\infty}\left(\int_\Omega\left\vert \nabla v_n\right\vert^2\,dx\right)^2=\liminf_{n\rightarrow+\infty}\mu_n=\mu_1.\nonumber
\end{equation}
On the other hand, $\int_\Omega\left( v_n\right)^4\,dx=1$ and $v_n\rightarrow v$ in $L^4(\Omega)$
imply that
$\int_\Omega v^4\,dx=1$. It follows that
\begin{equation}
\left(\int_\Omega\left\vert \nabla v\right\vert^2\,dx\right)^2 \leq\mu_1\int_\Omega v^4\,dx.\nonumber
\end{equation}
The above inequality and the variational characterization of $\mu_1$ imply that
\begin{equation}
\left(\int_\Omega\left\vert \nabla v\right\vert^2\,dx\right)^2 =\mu_1\int_\Omega v^4\,dx.\nonumber
\end{equation}
Then Proposition 3.2 follows that $v$ is positive or negative. Without loss of generality,
we may assume that $v>0$ in $\Omega$. Since $v_n\rightharpoonup v$ $X$,
going if necessary to a subsequence, we can assume that
\begin{equation}
v_n\rightarrow v\,\,\,\text{in}\,\, L^q(\Omega)\,\, \text{with}\,\, q\in\left(1,2^*\right),\nonumber
\end{equation}
\begin{equation}
v_n\rightarrow v\,\,\,\text{in a.e.}\,\,  \Omega.\nonumber
\end{equation}
So we conclude that
\begin{equation}
\left\vert \Omega_n^{-}\right\vert\rightarrow 0,\nonumber
\end{equation}
where $\Omega_n^{-}$ denotes the negative set of $u_n$.
This contradicts estimate (\ref{15}).\qed

\section{Positive solutions}

\quad\, In this section, we apply Theorem 1.1 and 1.2 to study the existence of positive solutions for (\ref{fp}).
\\ \\
\textbf{Lemma 4.1.} \emph{Assume (f1)--(f3) hold. Then $\left(1/f_0,0\right)$ is a bifurcation
point of (\ref{fp}) and the associated bifurcation branch $\mathcal{C}$ in $\mathbb{R}\times X$
whose closure contains $\left(1/f_0,0\right)$ is either unbounded or contains a pair $\left(\frac{\overline{\lambda}}{f_0\lambda_1}, 0\right)$ where $\overline{\lambda}$
is an eigenvalue of (\ref{le}) and $\overline{\lambda}\neq\lambda_1$.}
\\ \\
\textbf{Proof.}
Let $\vartheta:\Omega\times\mathbb{R}\setminus\mathbb{R}^-\rightarrow \mathbb{R}\setminus\mathbb{R}^-$ be a
Carath\'{e}odory function such that
\begin{equation}
f(x,s)=a\lambda_1f_0s+\vartheta(x,s)\nonumber
\end{equation}
with
\begin{equation}\label{ehac}
\lim_{s\rightarrow0^+}\frac{\vartheta(x,s)}{a\lambda_1s}=0\,\,\text{and}\,\,\lim_{ s\rightarrow+\infty}\frac{\vartheta(x,s)}{s^3}=b\mu_1f_\infty\,\,\text{uniformly with respect to a.e.\,\,} x\in\Omega.
\end{equation}
From (\ref{ehac}), we can see that $\lambda\vartheta$ satisfies the assumptions of (\ref{a1}) and (G). Now,
Theorem 1.1 can be applied to get the results of this lemma.
\qed\\

Let $P=\{u\in X:u(x)>0 \,\,\text{for all}\,\, x\in \Omega\}$ be the positive cone in $X$.
\\ \\
\textbf{Lemma 4.2.} \emph{We have $\mathcal{C}\subseteq \left(P\cup\left\{\left(1/f_0,0\right)\right\}\right)$ and the last alternative of Lemma 4.1 is
impossible.}
\\ \\
\textbf{Proof.} Firstly, by the strong maximum principle [\ref{GT}, Theorem 8.19], we know that $u>0$ in the whole domain for any nontrivial solution $(\lambda,u)\in \mathcal{C}$. So we have $\mathcal{C}\subseteq\left(P\cup\left(\mathbb{R}\times \{0\}\right)\right)$. Suppose on the contrary, if there exists $\left(\lambda_m,u_m\right)\rightarrow\left(\frac{\overline{\lambda}}{f_0\lambda_1},0\right)$
when $m\rightarrow+\infty$ with $\left(\lambda_m,u_m\right)\in \mathcal{C}$, $u_m \not\equiv 0$ and $\overline{\lambda}\neq\lambda_1$.
So we have $u_m\in P$ for each $m\in \mathbb{N}$.
Let $v_m :=u_m/\left\Vert u_m\right\Vert$, then $\left(\lambda_m,v_m\right)$ satisfies
\begin{equation}\label{evm}
v_m=G\left(\frac{\lambda_m}{a+b\left\Vert u_m\right\Vert^2}\left( a\lambda_1f_0v_m+\frac{\vartheta\left(\lambda,u_m\right)}{\left\Vert u_m\right\Vert}\right)\right).\nonumber
\end{equation}
By an argument similar to that of Theorem 1.1, we obtain that for some convenient
subsequence $v_m\rightarrow v_0$ as $m\rightarrow+\infty$. Now $v_0$ verifies the equation
\begin{equation}
-\Delta v= \overline{\lambda}v\nonumber
\end{equation}
and $\left\Vert v_0\right\Vert = 1$.
Hence $v_0$ must change its sign, and this is a contradiction. Furthermore, it follows that $\mathcal{C}\subseteq \left(P\cup\left\{\left(1/f_0,0\right)\right\}\right)$
and $\mathcal{C}$ is unbounded in $\mathbb{R}\times X$.\qed\\ \\
\textbf{Remark 4.1.} From the proof of Lemma 4.2, we can see that $\left(1/f_0, 0\right)$ is the unique bifurcation point of the set of positive solution of (\ref{fp}).\\

Next, we give a Sturm type comparison theorem.
\\ \\
\textbf{Theorem 4.1.} \emph{Assume that $g$ and $f_n$ be two weight functions with $\underset{n\rightarrow +\infty}{\lim}f_n(x)=+\infty$ for a.e. $x\in\Omega$, and satisfy
$g\in L^3(\Omega)$, $f_n\in L^3(\Omega)$ and $f_n\not\equiv g$ a.e. in $\Omega$ for any $n$ large enough.
Let $u$ be a positive weak solution of
\begin{equation}
\left\{
\begin{array}{l}
-\Vert u\Vert^2\Delta u=g(x)u^3\,\, \text{in}\,\,\Omega,\\
u=0~~~~~~~~~~~~~~~~~\,\text{on}\,\,\partial\Omega.
\end{array}
\right.\nonumber
\end{equation}
Then any solution $v_n\in W^{2,1}(\Omega)$ of
\begin{equation}
-\Vert v\Vert^2\Delta v=f_n(x)v^3\,\, \text{in}\,\,\Omega\nonumber
\end{equation}
must change sign for $n$ large enough.}
\\ \\
\textbf{Proof.} Without loss of generality, we assume that $\Vert u\Vert=\Vert v_n\Vert=1$ for any $n\in \mathbb{N}$. Suppose the contrary, we may assume that $v_n>0$ for $n$ large enough. Then by the following Picone's identity
\begin{equation}
\vert \nabla u\vert^2-\nabla \left(\frac{u^2}{v}\right)\nabla v=\vert \nabla u\vert^2+\frac{u^2}{v^2}\vert \nabla v\vert^2-2\frac{u}{v}\nabla u\nabla v\geq0\nonumber
\end{equation}
and an easy calculation, we obtain that
\begin{equation}
0\leq\int_\Omega\left(\vert \nabla u\vert^2+\frac{u^2}{v_n^2}\left\vert \nabla v_n\right\vert^2-2\frac{u}{v_n}\nabla u\nabla v_n\right)\,dx=\int_\Omega \left(g(x)u^2-f_n(x)v_n^2\right) u^2\,dx\nonumber
\end{equation}
for $n$ large enough.

For any $M>0$, we have that
\begin{equation}
1=\int_\Omega f_n(x) v_n^4\,dx\geq M\int_\Omega v_n^4\,dx\,\,\text{for $n$ large enough}.\nonumber
\end{equation}
It follows that $v_n\rightarrow 0$ in $L^4(\Omega)$. We claim that $v_n\leq u$ a.e. in $\Omega$ for $n$ large enough.
Otherwise, there exists $\Omega_0\subseteq \Omega$ with $\left\vert \Omega_0\right\vert>0$ such that $v_n>u$ in $\Omega_0$.
Then we have that
\begin{equation}
\int_{\Omega_0}u^4\,dx\leq \int_{\Omega_0}v_n^4\,dx\leq \int_{\Omega}v_n^4\,dx\rightarrow 0\,\, \text{as}\,\, n\rightarrow+\infty.\nonumber
\end{equation}
This is a contradiction.%$ f_n(x)\leq M_n$ such that $M_n\rightarrow\infty$
 Thus we have that
\begin{equation}
\int_\Omega \left(g(x)u^2-f_n(x)v_n^2\right) u^2\,dx=\int_\Omega f_n(x)v_n^2\left(v_n^2-u^2\right)\,dx\leq 0\nonumber
\end{equation}
for $n$ large enough.

Consequently we haver $u = cv_n$ for $n$ large enough. Furthermore, we have $u=\pm v_n$ for $n$ large enough. But this is impossible since $f_n\not\equiv g$ a.e. in $\Omega$ for $n$ large enough. This accomplishes the
proof.\qed
\\ \\
\textbf{Proof of Theorem 1.3.} In view of Lemma 4.1, 4.2 and Remark 4.1, it is sufficient to show that $\mathcal{C}$ joins $\left(1/f_0,0\right)$ to $\left(1/f_\infty,\infty\right)$.
Let $\left(\lambda_n, u_n\right) \in \mathcal{C}$ where $u_n\not\equiv 0$ satisfies $\lambda_n+\left\Vert u_n\right\Vert\rightarrow+\infty.$
Since (0,0) is the only solution of (\ref{fp}) for $\lambda = 0$, we have $\mathcal{C}\cap\left(\{0\}\times X\right)=\emptyset$. It follows that $\lambda_n >0$ for all $n \in \mathbb{N}$.

We claim that there exists a constant $M$ such that $\lambda_n
\in(0,M]$ for $n\in \mathbb{N}$ large enough. On the contrary, we suppose that $\lim_{n\rightarrow +\infty}\lambda_n=+\infty.$
Since  $\left(\lambda_n, u_n\right) \in \mathcal{C}$, it follows that
\begin{equation}
\Delta u_n+\frac{\lambda_n}{a+b\left\Vert u_n\right\Vert^2}\frac{f\left(x,u_n\right)}{u_n}u_n=0\,\, \text{in}\,\,\Omega.\nonumber
\end{equation}
We divide two cases to deduce a contradiction.

\emph{Case 1}. There exists a constant $c>0$ such that $\left\Vert u_n\right\Vert\leq c$ for $n$ large enough.

In this case, we have $\frac{1}{a+b\left\Vert u_n\right\Vert^2}\geq \frac{1}{a+bc^2}$. From (f1)--(f3), we can see that $\frac{f\left(x,u_n\right)}{u_n}\geq \sigma$ for some $\sigma>0$ and a.e. $x\in \Omega$ and all $n\in \mathbb{N}$.
By (f2) and (f3), we can obtain $F\left(x,u_n\right)/u_n\in L^{3/2}(\Omega)$, where $F$ denotes the usual Nemitsky operator associated with $f$.
By applying Theorem 2.6 of [\ref{AH}] with $N_0=3$, we have that $u_n$ must change its sign in $\Omega$, which contradicts Lemma 4.2.

\emph{Case 2}. $\left\Vert u_n\right\Vert\rightarrow +\infty$ as $n\rightarrow +\infty$.

Now, we consider
\begin{equation}
\left\Vert u_n\right\Vert^2\Delta u_n+\lambda_n\frac{\left\Vert u_n\right\Vert^2}{a+b\left\Vert u_n\right\Vert^2}\frac{f\left(x,u_n\right)}{u_n^3}u_n^3=0\,\, \text{in}\,\,\Omega.\nonumber
\end{equation}
For any fixed $\varepsilon\in\left(0,1/b\right)$, obviously, there exists $N_1>0$ such that
\begin{equation}\label{b3}
\frac{\left\Vert u_n\right\Vert^2}{a+b\left\Vert u_n\right\Vert^2}\geq \frac{1}{b}-\varepsilon
\end{equation}
for any $n>N_1$. By (f1)--(f3), there exists a constant $\varrho>0$ such that $\frac{f\left(x,u_n\right)}{u_n^3}\geq \varrho$ for a.e. $x\in \Omega$ and $n$ large enough.
Let
\begin{equation}
f_n(x)=\lambda_n\frac{\left\Vert u_n\right\Vert^2}{a+b\left\Vert u_n\right\Vert^2}\frac{f\left(x,u_n(x)\right)}{u_n(x)}\frac{1}{u_n^2(x)}.\nonumber
\end{equation}
Then we have $\lim_{n\rightarrow +\infty}f_n(x)=+\infty$ for a.e. $x\in\Omega$. By (f1)--(f3) and [\ref{GT}, Theorem 8.12], we know that
$u_n\in W^{2,2}(\Omega)$. So we have $u_n\in W^{2,1}(\Omega)$.

On the other hand, by Theorem 1.1 of [\ref{FZZ}], we have that for any nonempty compact subset $K\subseteq \Omega$, there exists a positive constant $c_n$
such that $u_n\geq c_n$ a.e. in $K$. Then it is not difficult to show that $f_n\in L^3(K)$ for any fixed $n\in \mathbb{N}$.
Applying Theorem 4.1 on $K$ with $g(x)\equiv \mu_1$, we have that $u_n$ must change its sign in $K$ for $n$ large enough. This is a contradiction.

Therefore, we get that
\begin{equation}
\left\Vert u_n\right\Vert\rightarrow+\infty.\nonumber
\end{equation}
Let $\varpi:\Omega\times\mathbb{R}\setminus\mathbb{R}^-\rightarrow \mathbb{R}\setminus\mathbb{R}^-$ be a
Carath\'{e}odory function such that
\begin{equation}
f(x,s)=b\mu_1f_\infty s^3+\varpi(x,s)\nonumber
\end{equation}
with
\begin{equation}\label{ehac2}
\lim_{s\rightarrow+\infty}\frac{\varpi(x,s)}{s^3}=0\,\,\text{and}\,\,\lim_{ s\rightarrow0^+}\frac{\varpi(x,s)}{s}=a\lambda_1 f_0\,\,\text{uniformly with respect to a.e.\,\,} x\in\Omega.
\end{equation}
Then $\left(\lambda_n,u_n\right)$ satisfies
\begin{equation}\label{evm2}
u_n=G\left(\frac{\lambda_n}{a+b\left\Vert u_n\right\Vert^2}\left( b\mu_1f_\infty u_n^3+\varpi\left(x,u_n\right)\right)\right).\nonumber
\end{equation}
Dividing the above equation by $\left\Vert u_n\right\Vert$ and letting $\overline{u}_n=u_n/\left\Vert u_n\right\Vert$, we get that
\begin{equation}
\overline{u}_n=G\left(\frac{\lambda_n \left\Vert u_n\right\Vert^2}{a+b\left\Vert u_n\right\Vert^2}\left( b\mu_1f_\infty \overline{u}_n^3+\frac{\varpi\left(x,u_n\right)}{\left\Vert u_n\right\Vert^3}\right)\right).\nonumber
\end{equation}
\indent Next, we show that
\begin{equation}\label{44}
\lim_{n\rightarrow+\infty}\frac{ \varpi\left(x,u_n(x)\right)}{\left\Vert u_n\right\Vert^{3}}=0\,\, \text{in\,\,} L^{q'}(\Omega)
\end{equation}
for some $q<2^*$. Without loss of generality, we may assume that $q\geq4$. Otherwise, we can consider $\tilde{q}=cq$, $c>1$ such that $\tilde{q}\in\left[4,2^*\right)$.

From (\ref{ehac2}), for any $\varepsilon>0$, we can choose $\delta=\delta(\varepsilon)$ and $M=M(\delta)$ such that for a.e. $x\in \Omega$ and any $n\in \mathbb{N}$, the following relation hold:
\begin{equation}\label{ehC}
\vert\varpi\left(x,s\right)\vert\leq M\,\,\,\text{for}\,\, \vert s\vert\leq \delta
\end{equation}
and
\begin{equation}\label{ehC11}
\vert\varpi\left(x,s\right)\vert\leq \epsilon \vert s\vert^{3}\,\,\,\text{for}\,\, \vert s\vert> \delta.
\end{equation}
By (\ref{ehC}) and (\ref{ehC11}), we can easily show that
\begin{equation}
\int_\Omega \left\vert \frac{ \varpi\left(x,u_n(x)\right)}{\left\Vert u_n\right\Vert^{3}}\right\vert\,dx\leq \frac{c}{\left\Vert u_n\right\Vert^{3q'}}+\varepsilon\int_\Omega \left\vert \overline{u}_n\right\vert^{3q'}\,dx.\nonumber
\end{equation}
It follows from $q\geq 4$ that $3q'<2^*$.
From the above inequality, $u_n\rightarrow+\infty$ in $X$, $\left\Vert \overline{u}_n\right\Vert=1$, we can get the desired result.
It is not difficult to show that $\overline{u}_n^3$ is bounded in $L^{q'}(\Omega)$ and $\frac{\lambda_n \left\Vert u_n\right\Vert^2}{a+b\left\Vert u_n\right\Vert^2}$ is bounded in $\mathbb{R}$ for $n$ large enough.
By the compactness of $G$ we obtain that
\begin{equation}
-\left\Vert \overline{u}\right\Vert^2\Delta \overline{u}=\overline{\mu}\mu_1f_\infty \overline{u}^3,\nonumber
\end{equation}
where $\overline{u}=\underset{n\rightarrow+\infty}\lim \overline{u}_n$ and $\overline{\mu}=\underset{n\rightarrow+\infty}\lim\lambda_n$, again choosing a subsequence and relabeling it if necessary.

It is clear that $\left\Vert \overline{u}\right\Vert=1$ and $\overline{u}\in \overline{\mathcal{C}}\subseteq \mathcal{C}$ since $\mathcal{C}$
is closed in $\mathbb{R}\times X$.
Therefore, Lemma 4.2 implies $\overline{u}>0$. Theorem 1.2 shows that $\overline{\mu}=1/f_\infty$.
Therefore, $\mathcal{C}$ joins $\left(1/f_0,0\right)$ to $\left(1/f_\infty,\infty\right)$.
\qed
\\ \\
\textbf{Corollary 4.1.} \emph{Assume that $N\leq 3$ and $f$ satisfies (f1)--(f2). Then for
\begin{equation}
\lambda\in\left(1/f_0,1/f_\infty\right)\cup \left(1/f_\infty,1/f_0\right),\nonumber
\end{equation}
problem (\ref{fp}) possesses at least one positive solution. In particularly, if $f_0>1$ ($<1$) and $f_\infty<1$ ($>1$) then (\ref{fp}) possesses at least one positive solution with $\lambda=1$.}

\section{Answer to an open problem}

\quad\, In this section, we consider problem (\ref{fp}) with $a=0$, i.e., the following problem
\begin{equation}\label{a01}
\left\{
\begin{array}{l}
-b\int_\Omega \vert \nabla u\vert^2\,dx\Delta u=\lambda f(x,u)\,\,\text{in}\,\, \Omega,\\
u=0~~~~~~~~~~~~~~~~~~~~~~~~~~~~~~~\,\text{on}\,\,\Omega.
\end{array}
\right.
\end{equation}
To prove Theorem 1.4, we recall the following topological lemma (see [\ref{MA}]).
\\ \\
\textbf{Lemma 5.1.} \emph{Let $\mathcal{X}$ be a Banach space with normal $\Vert\cdot\Vert_{\mathcal{X}}$ and let $\mathcal{C}_n$ be a family of closed connected subsets of $\mathcal{X}$. Assume that:}
\\

(i) \emph{there exist $z_n\in \mathcal{C}_n$, $n=1,2,\ldots$, and $z^*\in \mathcal{X}$, such that $z_n\rightarrow z^*$;}

(ii) \emph{$r_n=\sup \left\{\Vert x\Vert_{\mathcal{X}}\big| x\in \mathcal{C}_n\right\}=+\infty$;}

(iii) \emph{for every $R>0$, $\left(\cup_{n=1}^{+\infty} \mathcal{C}_n\right)\cap B_R$ is a relatively compact set of $\mathcal{X}$, where}
\begin{equation}
B_R=\left\{x\in \mathcal{X}|\Vert x\Vert_{\mathcal{X}}\leq R\right\}.\nonumber
\end{equation}
\emph{Then there exists an unbounded component $\mathcal{C}$ of $\mathfrak{D} =: \limsup_{n\rightarrow +\infty}\mathcal{C}_n$ and $z^*\in \mathcal{C}$.}
\\

In order to apply Lemma 5.1 to prove Theorem 1.4, we need to
discuss the following Kirchhoff-Laplace operator
\begin{equation}
\Phi_k(u):=-b\Vert u\Vert^2\Delta u.\nonumber
\end{equation}
Denote
\begin{equation}
\Phi(u):=\frac{1}{4}b\Vert u\Vert^4.\nonumber
\end{equation}
It is obvious
that the functional $\Phi$ is continuously G\^{a}teaux differentiable whose G\^{a}teaux derivative at the point $u\in X$ is the functional
$\Phi'(u)\in X^*$, given by
\begin{equation}
\left\langle\Phi'(u),v\right\rangle=b\Vert u\Vert^2\int_\Omega \nabla u\cdot\nabla v\,dx.\nonumber
\end{equation}
Obviously, the Kirchhoff-Laplace operator is the derivative operator
of $\Phi$ in the weak sense. We have the following properties about the derivative operator of $\Phi$.
\\ \\
\textbf{Proposition 5.1.} \emph{Let $L=\Phi'$. Then we have}\\

(i) \emph{$L: X\rightarrow X^*$ is a continuous and strictly monotone operator;}

(ii) \emph{$L$ is a mapping of type $\left(S_+\right)$, i.e., if $u_n\rightharpoonup u$ in $X$ and $\underset{n\rightarrow+\infty}{\overline{\lim}}
\left\langle L\left(u_n\right)-L(u),u_n-u\right\rangle\leq0$, then $u_n \rightarrow u$ in $X$;}

(iii) \emph{$L(u):X\rightarrow X^*$ is a homeomorphism.}
\\
\\
\textbf{Proof.}
(i) Let $u_n\rightarrow u$ in $X$, i.e. $\lim_{n\rightarrow+\infty}\left\Vert u_n-u\right\Vert=0$. Then we can easily see that
$\lim_{n\rightarrow+\infty}\left\Vert u_n\right\Vert=\left\Vert u\right\Vert$. For any $v\in X$, by using of the H\"{o}lder's inequality, we have that
\begin{eqnarray}
\left\vert\left\langle L\left(u_n\right)-L(u),v\right\rangle \right\vert&=&\left\vert b\left\Vert u_n\right\Vert^2\int_\Omega \nabla u_n\cdot\nabla v\,dx-b\Vert u\Vert^2\int_\Omega \nabla u\cdot\nabla v\,dx\right\vert\nonumber \\
&=&\left\vert b\left\Vert u_n\right\Vert^2\int_\Omega \left(\nabla u_n-\nabla u\right)\cdot\nabla v\,dx+b\left(\left\Vert u_n\right\Vert^2-\left\Vert u\right\Vert^2\right)\int_\Omega \nabla u\cdot\nabla v\,dx\right\vert\nonumber\\
&\leq& b\left\Vert u_n\right\Vert^2\left\Vert u_n-u\right\Vert \Vert v\Vert+b\left\vert\left\Vert u_n\right\Vert^2-\left\Vert u\right\Vert^2\right\vert \Vert u\Vert\Vert v\Vert\nonumber\\
&\rightarrow &0\,\,\text{as}\,\, n\rightarrow+\infty.\nonumber
\end{eqnarray}
It follows that $L$ is a continuous operator.

For any $u$, $v\in X$ with $u\neq v$ in $X$, by the Cauchy's inequality, we obtain that
\begin{eqnarray}\label{au1}
\left\langle L\left(u\right)-L(v),u-v\right\rangle&=&\langle L(u),u\rangle-\langle L(u),v\rangle-\langle L(v),u\rangle+\langle L(v),v\rangle\nonumber\\
&=& b\Vert u\Vert^4-b \Vert v\Vert^2\int_\Omega\nabla u\cdot\nabla v\,dx-b \Vert u\Vert^2\int_\Omega\nabla u\cdot\nabla v\,dx+b\Vert v\Vert^4\nonumber\\
&=&b\left(\Vert u\Vert^4+\Vert v\Vert^4- \left(\Vert u\Vert^2+\Vert v\Vert^2\right)\int_\Omega\nabla u\cdot\nabla v\,dx\right)\nonumber\nonumber\\
&\geq&b\left(\Vert u\Vert^4+\Vert v\Vert^4- \frac{\left(\Vert u\Vert^2+\Vert v\Vert^2\right)^2}{2}\right)\nonumber\\
&=&b\left(\Vert u\Vert^2-\Vert v\Vert^2\right)^2\nonumber\\
&\geq& 0,
\end{eqnarray}
i.e. $L$ is monotone. In fact $L$ is strictly monotone. Indeed, if $\langle L(u)-L(v),u-v\rangle=0$, then from (\ref{au1}) we have that
\begin{equation}
\Vert u\Vert=\Vert v\Vert,\,\,  \nabla u\equiv \nabla v\,\,\text{a.e. in}\,\, \Omega.\nonumber
\end{equation}
It follows that $\Vert u-v\Vert=0$, i.e. $u\equiv v$,
which is contrary with $u\neq v$  in $X$.
Therefore, $\langle L(u)-L(v),u-v\rangle>0$. It follows that $L$ is a strictly monotone operator in $X$.

(ii) From (i), if $u_n\rightharpoonup u$ and $\underset{n\rightarrow+\infty}{\overline{\lim}}
\left\langle L\left(u_n\right)-L(u),u_n-u\right\rangle\leq0$, then we have that
\begin{equation}
\underset{n\rightarrow+\infty}{\lim}
\left\langle L\left(u_n\right)-L(u),u_n-u\right\rangle=0.
\nonumber
\end{equation}
In view of (\ref{au1}), we obtain that $\nabla u_n$ converges in measure to $\nabla u$ in $\Omega$, so
we get a subsequence (which we still denote by $u_n$) satisfying $\nabla u_n(x)\rightarrow \nabla u (x)$, a.e. $x\in \Omega$.
Moreover, we also have that 	
\begin{equation}\label{au2}
\underset{n\rightarrow+\infty}{\lim}\int_\Omega
\left\vert \nabla u_n\right\vert ^{2}\,{d}x=\int_\Omega
\left\vert \nabla u\right\vert ^{2}\,{d}x.
\end{equation}
From (\ref{au2}) it follows that the integrals of the functions family $\left\{\left\vert \nabla u_n\right\vert^2\right\}$
possess absolutely equi-continuity on $\Omega$ (see [\ref{N}, Chapter 6, Section 3]).
Since
\begin{equation}
\left\vert \nabla u_n-\nabla u\right\vert^2\leq 2\left(\left\vert \nabla u_n\right\vert^2+\vert \nabla u\vert^2\right),\nonumber
\end{equation}
the integrals of the family $\left\{\left\vert \nabla u_n-\nabla u\right\vert^2\right\}$
are also absolutely equi-continuous
on $\Omega$ and therefore we have
\begin{equation}
\underset{n\rightarrow+\infty}{\lim}\int_\Omega \left\vert \nabla u_n-\nabla u\right\vert ^{2}\,{d}x=0.\nonumber
\end{equation}
Therefore, $u_n\rightarrow u$, i.e. $L$ is of type $\left(S_+\right)$.

(iii) It is clear that $L$ is an injection. Since
\begin{equation}
\lim_{\Vert u\Vert\rightarrow+\infty}\frac{\langle L(u),u\rangle}{\Vert u\Vert}
=\lim_{\Vert u\Vert\rightarrow+\infty}b\Vert u\Vert^3=+\infty,\nonumber
\end{equation}
$L$ is coercive, thus $L$ is a surjection in view of Minty-Browder Theorem (see [\ref{Z},
Theorem 26A]). Hence $L$ has an inverse mapping $L^{-1}:X^*\rightarrow X$. Therefore, the continuity
of $L^{-1}$ is sufficient to ensure $L$ to be a homeomorphism.

If $f_n$, $f\in X^*$, $f_n\rightarrow f$, let $u_n =L^{-1}\left(f_n\right)$, $u = L^{-1}(f)$,
then $L\left(u_n\right) = f_n$, $L(u) = f$.
The coercive property of $L$ implies that $\left\{u_n\right\}$ is bounded in $X$.
We can assume that $u_{n_k}\rightharpoonup u_0$ in $X$.
By $f_{n_k} \rightarrow f$ in $X^*$, we have
\begin{equation}
\lim_{k\rightarrow+\infty}\left\langle L\left(u_{n_k}\right)-L\left(u_0\right),u_{n_k}-u_0\right\rangle
=\lim_{n\rightarrow+\infty}\left\langle f_{n_k},u_{n_k}-u_0\right\rangle=0.\nonumber
\end{equation}
Since $L$ is of type $\left(S_+\right)$, $u_{n_k}\rightarrow u_0$.
Furthermore, the continuity of $L$ implies that $L\left(u_0\right)=L(u)$.
By injectivity of $L$, we have $u_0 = u$. So $u_{n_k}\rightarrow u$.
We claim that $u_n\rightarrow u$ in $X$. Otherwise, there would exist a
subsequence $\left\{u_{m_j}\right\}$ of $\left\{u_n\right\}$ in $X$ and an $\varepsilon_0 > 0$, such that for any $j\in \mathbb{N}$,
we have $\left\Vert u_{m_j}-u\right\Vert\geq \varepsilon_0$. But reasoning as above, $\left\{u_{m_j}\right\}$ would contain a
further subsequence $u_{m_{j_l}}\rightarrow u$ in $X$ as $l \rightarrow+\infty$, which is a contradiction to
$\left\Vert u_{m_{j_l}}-u\right\Vert\geq \varepsilon_0$. Therefore, $L^{-1}$ is continuous.\qed
\\ \\
\textbf{Remark 5.1.} Note that in the proof of Lemma 4.3 of [\ref{LLS}], the authors deduced $\left\langle u_n,u_n\right\rangle\rightarrow \langle u,u\rangle$
from $\left(a+b\left\Vert u_n\right\Vert^2\right)\left[\left\langle u_n,u_n\right\rangle-\left\langle u_n,u\right\rangle\right]\rightarrow 0$. Clearly, this is right for the case of $a>0$;
but it may be not right for the case of $a=0$ because $\left\Vert u_n\right\Vert$ may converge to 0. However, by Lemma 5.2 (ii) we can get $\left\langle u_n,u_n\right\rangle\rightarrow \langle u,u\rangle$ from $b\left\Vert u_n\right\Vert^2\left[\left\langle u_n,u_n\right\rangle-\left\langle u_n,u\right\rangle\right]\rightarrow 0$ immediately.
So we complete the proof of Lemma 4.3 of [\ref{LLS}] at here.
\\

Now, consider the following auxiliary problem
\begin{equation}\label{ap2}
\left\{
\begin{array}{l}
-b\Vert u\Vert^2\Delta u=f(x)\,\, \text{in}\,\,\Omega,\\
u=0~~~~~~~~~~~~~~~~\,\,\,\text{on}\,\,\partial\Omega.
\end{array}
\right.
\end{equation}
\noindent\textbf{Lemma 5.2.} For any $f\in X^*$, problem (\ref{ap2}) has a unique weak solution.
\\ \\
\textbf{Proof.} Define $\mathcal{F}(v)=\int_\Omega fv\,dx$ for any $v\in X$. We can easily see $\mathcal{F}$ is a continuous linear functional on
$X$. Since $L$ is a homeomorphism, (\ref{ap2}) has a unique
solution.\qed\\

Let us denote by $S(f)$ the unique weak solution of (\ref{ap2}). Then $S:X^*\rightarrow X$ is a continuous operator. Since the embedding of $X\hookrightarrow L^q(\Omega)$ is compact for each $q\in\left(1,2^*\right)$,
the restriction of $S$ to $L^{q'}(\Omega)$ is a completely operator.

For any $n\in \mathbb{N}$, we study the following auxiliary problem
\begin{equation}\label{an1}
\left\{
\begin{array}{l}
-\left(\frac{1}{n}+b\int_\Omega \vert \nabla u\vert^2\,dx\right)\Delta u=\lambda f(x,u)\,\,\text{in}\,\, \Omega,\\
u=0~~~~~~~~~~~~~~~~~~~~~~~~~~~~~~~~~~~~~~~~~\text{on}\,\,\Omega.
\end{array}
\right.
\end{equation}
Clearly, we can see that $f_0=\widetilde{f}_0 n$.
Theorem 1.3 implies
that there exists a sequence of unbounded continua $\mathcal{C}_n$ of positive solutions to problem (\ref{an1})
emanating from $\left(\frac{1}{\widetilde{f}_0 n}, 0\right)$
and joining to $\left(\frac{1}{f_\infty}, +\infty\right)$.
\\ \\
\textbf{Proof of Theorem 1.4.} Let ${\mathcal{X}}=\mathbb{R}\times X$ under the
product topology. Clearly, ${\mathcal{X}}$ is a Banach space.
Now, we verify the assumptions of Lemma 5.1. Taking $z_n=\left(\frac{1}{\widetilde{f}_0 n}, 0\right)$ and $z^*=\left(0, 0\right)$,
we have that $z_n\rightarrow z^*$. So (i) is satisfied. (ii) is obvious.

By (f3) and (f4), we can show that $F(x,u(x))\in L^{4/3}(\Omega)$. Now, we divide two possibilities to verify (iii): (a) $n<+\infty$ and (b) $n=+\infty$. If the case (a) occurs, the
completely continuous property of $G$ implies (iii). If the case (b) occurs, the compactness of $S$ follows (iii).
By Lemma 5.1, there exists an unbounded component $\mathcal{C}$ of $\lim\sup_{n\rightarrow+\infty} \mathcal{C}_n$
such that $\left(0, 0\right)\in\mathcal{C}$ and $\left(1/f_\infty, +\infty\right)\in\mathcal{C}$.
This completes the proof.\qed
\\ \\
\textbf{Corollary 5.1.} \emph{Assume that $N\leq 3$ and $f$ satisfies (f1), (f3) and (f4). Then for
\begin{equation}
\lambda\in\left(0,1/f_\infty\right),\nonumber
\end{equation}
problem (\ref{a01}) possesses at least one positive solution. In particularly, if $f_\infty<1$ then (\ref{a01}) possesses at least one positive solution with $\lambda=1$.}
\\ \\
\textbf{Remark 5.2.} From Remark 1.4, we know that any weak solution of (\ref{an1}) belongs to $C^{1,\alpha}(\overline{\Omega})$ with $\alpha\in(0,1)$ under the assumptions of (f1), (f3) and (f4). It follows that $u\in C^{1,\alpha}(\overline{\Omega})$ with $\alpha\in(0,1)$ for any $(\lambda,u)\in \mathcal{C}$ which is obtained in Theorem 1.4.

\section{Conclusions and future work}

\quad\, This paper performs studies on the global bifurcation phenomena for the Kirchhoff type equations and its applications.
We firstly use bifurcation method to study the existence of positive solutions for the Kirchhoff type equations.
We give a positive answer to an open problem. Moreover, we also sharply obtain some important properties of the first eigenvalue of a nonlocal problem.

For future work, we plan to: 1) study the case of $f_0\not\in(0,+\infty)$ $(\widetilde{f}_0\not\in(0,+\infty))$ or $f_\infty\not\in(0,+\infty)$.
This plan comes from Theorem 1.2 and the following special example
\begin{equation}\label{61}
\left\{
\begin{array}{l}
-\left(a+b\int_\Omega \vert \nabla u\vert^2\,dx\right)\Delta u=\lambda u\,\,\text{in}\,\, \Omega,\\
u=0~~~~~~~~~~~~~~~~~~~~~~~~~~~~~~~~~\,\text{on}\,\,\Omega.
\end{array}
\right.
\end{equation}
It is easy to see that $f_\infty=0$. While, the positive solution pairs of (\ref{61}) must be
\begin{equation}
\left(\lambda_1\left(a+b\left\Vert \varphi_1\right\Vert^2\right),\varphi_1\right),\nonumber
\end{equation}
where
$\varphi_1$ is the corresponding positive eigenfunction to $\lambda_1$. Thus $\left(\lambda_1,0\right)$ is a bifurcation point of the set of positive solution of (\ref{61}) and $\left(\lambda_1\left(a+b\left\Vert \varphi_1\right\Vert^2\right),\varphi_1\right)$ is the corresponding unbounded branch.
2) Study the unilateral global bifurcation phenomenon and the existence of one-sign and sign-changing solutions for (\ref{bp}).


\begin{thebibliography}{99}
\bibitem{}\label{AH} W. Allegretto and Y.X. Huang, A Picone's identity for the $p$-Laplacian and applications,
Nonlinear Anal. 32 (1998), 819--830.

\bibitem{}\label{AP} A. Arosio and S. Pannizi, On the well-posedness of the Kirchhoff
string, Trans. Amer. Math. Soc. 348 (1996), 305--330.

\bibitem{}\label{APS} G. Autuori, P. Pucci and M. C. Salvatori, Asymptotic stability for anisotropic
Kirchhoff systems, J. Math. Anal. Appl. 352 (2009), 149--165.

\bibitem{}\label{CCS} M.M. Cavalcante, V.N. Cavalcante and J.A. Soriano, Global existence and
uniform decay rates for the Kirchhoff-Carrier equation with
nonlinear dissipation, Adv. Differential Equations 6 (2001),
701--730.

\bibitem{}\label{CKW} C. Chen, Y. Kuo and T. Wu, The Nehari manifold for a Kirchhoff type problem involving sign-changing weight
functions, J. Differential Equations 250 (4) (2011), 1876--1908.

\bibitem{}\label{CW} B. Cheng and X. Wu, Existence results of positive solutions of Kirchhoff type problems, Nonlinear Anal. 71 (10) (2009), 4883--4892.

\bibitem{}\label{CF} F.J.S.A. Corr\^{e}a and G.M. Figueiredo, On a elliptic
equation of $p$-kirchhoff type via variational methods, Bull.
Austral. Math. Soc. 74 (2006), 263--277.

\bibitem{}\label{DH} G. Dai and R. Hao, Existence of solutions for a $p(x)$-Kirchhoff-type equation,
J. Math. Anal. Appl. 359 (2009), 275--284.

\bibitem{}\label{DJ} G. Dai and J. Wei, Infinitely many non-negative solutions for a $p(x)$-Kirchhoff-type
 problem with Dirichlet boundary condition, Nonlinear Anal. 73 (2010), 3420--3430.

\bibitem{}\label{DL} G. Dai and D. Liu, Infinitely many positive solutions for a $p(x)$-Kirchhoff-type
equation, J. Math. Anal. Appl. 359 (2009), 704--710.

\bibitem{}\label{DM} G. Dai and R. Ma,
Solutions for a $p(x)$-Kirchhoff type equation with Neumann boundary data, Nonlinear Anal. Real World Appl. 12 (2011), 2666--2680.

\bibitem{}\label{DS} P. D'Ancona and S. Spagnolo, Global solvability for the degenerate
Kirchhoff equation with real analytic data, Invent. Math. 108 (1992),
247--262.

\bibitem{}\label{DS1} P. D'Ancona and Y. Shibata, On global solvability of nonlinear viscoelastic equations in the analytic category, Math. Methods Appl.
Sci. 17 (6) (1994), 477--486.

\bibitem{}\label{DPM} M. Del Pino and R. Man\'{a}sevich, Global bifurcation from the eigenvalues of the $p$-Lapiacian, J. Differential Equations 92 (1991), 226--251.

\bibitem{}\label{D} M. Dreher, The Kirchhoff equation for the $p$-Laplacian, Rend. Semin. Mat.
Univ. Politec. Torino 64 (2006), 217--238.

\bibitem{}\label{D1} M. Dreher, The ware equation for the $p$-Laplacian, Hokkaido Math. J. 36 (2007), 21--52.

\bibitem{}\label{E} L.C. Evans, Partial differential equations, AMS, Rhode Island, 1998.

\bibitem{}\label{F} X.L. Fan, On nonlocal $p(x)$-Laplacian Dirichlet problems, Nonlinear Anal.
72 (2010), 3314--3323.

\bibitem{}\label{F1} X.L. Fan, Global $C^{1,\alpha}$ regularity for variable exponent elliptic
equations in divergence form, J. Differential Equations 235 (2007), 397--417.

\bibitem{}\label{FZ} X.L. Fan and D. Zhao, A class of De Giorgi type and H\"{o}lder continuity, Nonlinear Anal. 36 (1996), 295--318.

\bibitem{}\label{FZZ} X.L. Fan, Y.Z. Zhao and Q.H. Zhang, A strong maximum principle for
$p(x)$-Laplace equations, Chinese J. Contemp. Math. 24 (3) (2003),
277--282.

\bibitem{}\label{GT} D. Gilbarg and N.S. Trudinger, Elliptic partial differential equations of second order, Springer, Berlin, 2001.

\bibitem{}\label{HZ} X. He and W. Zou, Infinitely many positive solutions for Kirchhoff-type problems,
Nonlinear Anal. 70 (2009), 1407--1414.

\bibitem{}\label{K} G. Kirchhoff, Mechanik, Teubner, Leipzig, 1883.

\bibitem{}\label{LLS} Z. Liang, F. Li and J. Shi, Positive solutions to Kirchhoff type equations with nonlinearity
having prescribed asymptotic behavior, Ann. I. H. Poncar\'{e}-AN (2013), http://dx.doi.org/10.1016/j.anihpc.2013.01.006.

\bibitem{}\label{L} J.L. Lions, On some equations in boundary value problems of mathematical
physics, in: Contemporary Developments in Continuum Mechanics and
Partial Differential Equations (Proc. Internat. Sympos., Inst. Mat.
Univ. Fed. Rio de Janeiro, Rio de Janeiro, 1977), in: North-Holland
Math. Stud., vol. 30, North-Holland, Amsterdam, 1978, pp. 284--346.

\bibitem{}\label{MA} R. Ma and Y. An, Global structure of positive solutions for nonlocal boundary value
problems involving integral conditions, Nonlinear Anal. 71 (2009, 4364--4376.

\bibitem{}\label{MR} T.F. Ma and J.E.Mu\~{n}oz Rivera, Positive solutions for a nonlinear nonlocal elliptic transmission problem, Appl. Math. Lett. 16 (2) (2003), 243--248.

\bibitem{}\label{MZ} A. Mao and Z. Zhang, Sign-changing and multiple solutions of Kirchhoff type problems without the P.S. condition, Nonlinear
Anal. 70 (3) (2009), 1275--1287.

\bibitem{}\label{M} C.B. Morry, Multiple integrals in the calculus of variations, Springer-Verlag, New York, 1966.

\bibitem{}\label{N} I.P. Natanson, Theory of Functions of a Real Variable, Nauka, Moscow, 1950.

\bibitem{}\label{PZ} K. Perera and Z. Zhang, Nontrivial solutions of Kirchhoff-type problems via the Yang index, J. Differential Equations 221 (1) (2006),
246--255.

\bibitem{}\label{R1} P.H. Rabinowitz, Some aspects of nonlinear eigenvalue problems,
Rocky Mountain J. Math. 3(2) (1973), 161--202.

\bibitem{}\label{R2} P.H. Rabinowitz, Some global results for nonlinear eigenvalue problems,
J. Funct. Anal. 7 (1971), 487--513.

\bibitem{}\label{ST} J. Sun and C. Tang, Existence and multiplicity of solutions for Kirchhoff type equations, Nonlinear Anal. 74 (4) (2011), 1212--1222.

\bibitem{}\label{Z} E. Zeidler, Nonlinear functional analysis and its applications,
Vol. II/B. Berlin-Heidelberg-New York 1985.

\bibitem{}\label{ZP} Z. Zhang and K. Perera, Sign changing solutions of Kirchhoff type problems via invariant sets of descent flow, J. Math. Anal.
Appl. 317 (2) (2006), 456--463.

\end{thebibliography}
\end{document}